\documentclass[12pt]{amsart}

\usepackage{tikz}  
\usetikzlibrary{decorations.markings}

\usepackage[colorlinks=true]{hyperref}
\usepackage{doi}

\newtheorem{theorem}{Theorem}

\begin{document}
\title{Cauchy's Residue Sore Thumb}
\author{Harold P. Boas}
\address{Department of Mathematics \\ Texas A\&M University \\ College Station TX 77843-3368 \\ USA}
\email{boas@tamu.edu}

\maketitle

Are you good at computing integrals? Try this one:
\begin{equation}
\int_{0}^{\infty} e^{\cos(x)} \sin(\sin(x)) \frac{x}{x^{2}+1} \,dx.
\label{eq:quiz}
\end{equation}
No fair peaking at the answer! But you get partial credit for showing at least that this improper integral converges.

If you find this problem a hard nut to crack, you are in good company. The integral is absent from the exhaustive tables~\cite{table} of Gradshteyn and Ryzhik,\footnote{Entry 3.973 is a near miss.} 
and when I fed this problem to Maple~18 and to \textit{Mathematica}~11, both software programs choked. Even the great Augustin-Louis Cauchy (1789--1857), who posed the problem, got the answer wrong on his first try. 

There is no hope to evaluate~\eqref{eq:quiz} by first computing an explicit antiderivative of the integrand. The failure of standard computer programs to produce an elementary antiderivative is compelling evidence that there is none; skeptical readers can prove the nonexistence by invoking the theory of differential fields as illustrated in the expository article~\cite{rosenlicht} by Maxwell Rosenlicht (1924--1999). 

Mathematicians of the 19th century knew so many special tricks for evaluating definite integrals that the Dutch scholar David Bierens de Haan (1822--1895) could write a book~\cite{expose} on the topic.
The approach to the integral~\eqref{eq:quiz} that likely occurs to a mathematician of the 21st century is Cauchy's residue theorem for functions of a complex variable. Indeed, the application of complex analysis to solve purely real problems is nowadays a familiar idea~\cite{laxzalcman}. 

Euler's formula for the complex exponential function implies that
\begin{equation*}
e^{e^{ix}}
=e^{\cos(x)+i\sin(x)} = e^{\cos(x)}\left( \cos(\sin(x)) + i \sin(\sin(x))\right).
\end{equation*}
Accordingly, a conceivable method for attacking~\eqref{eq:quiz} is first to integrate the expression \(e^{e^{iz}} z/(z^{2}+1)\) around a suitable  contour and then to take the imaginary part of the result. The integral over a simple, closed, counterclockwise curve in the upper half-plane surrounding the singular point~\(i\) equals \(2\pi i\) times the  residue---the coefficient of \( (z-i)^{-1}\) in the local expansion of the function in positive and negative powers of \( (z-i) \). The integrand can be written as \(  \left[ e^{e^{iz}} z(z+i)^{-1} \right] (z-i)^{-1}\), so the residue equals \(e^{e^{i^{2}}} i(2i)^{-1}\), whence the integral over the closed curve equals \(\pi i e^{1/e}\).
The nontrivial obstacle to executing this method is that the original integration path is not a closed curve. 

The integral~\eqref{eq:quiz} appears~\cite[appendix, formula~(25)]{Cauchy1825} in a long list of integrals that Cauchy evaluates\footnote{I write about Cauchy's work in the present tense, on the grounds that his mathematics is timeless.}  
in a memoir of 1825 that is often viewed as the origin of the residue theorem. Observing by symmetry that \eqref{eq:quiz}~is half the integral over the whole real line, Cauchy deduces the incorrect value \(\frac{1}{2} \pi e^{1/e}\) for~\eqref{eq:quiz}. He soon corrects the mistake \cite[p.~139]{correction}, acknowledging that a mishandling of his powerful new tool is the cause of his smashed thumb.\footnote{Cauchy's confession of error appears in a footnote. The hammer metaphor is my own.} 

The influential philosopher Imre Lakatos (1922--1974) emphasizes in a seminal book \cite{lakatos} that new concepts and theorems are generated by proofs---not the other way around. Two examples that he cites from the work of Cauchy are a purported proof of Euler's formula for polyhedra and a purported proof of continuity of convergent series of continuous functions, both of which Lakatos views as good arguments in search of valid theorems. Similarly, the residue calculus represents a remarkably successful technique, even though Cauchy's implementations lack accurate hypotheses. 

My goal in this article is not merely to supply a sound calculation of the integral~\eqref{eq:quiz} but also to formulate and prove natural theorems that realize Cauchy's original vision. My treatment differs from the standard exposition of the residue calculus in modern sources. 
So many cooks have seasoned the residue broth during the past two centuries that the recipe now has become codified in a form that loses sight of Cauchy's simple initial conception. 
I present in a few paragraphs\footnote{Amusingly, a catchphrase of the prolific Cauchy is ``en peu de mots'' (literally, ``in a few words'').} 
a self-contained development of the part of Cauchy's theory needed for evaluating integrals over the real axis. 

Although I do not aim to compete with comprehensive studies of the history and applications of Cauchy's work on complex integration (such as \cite{mitrinovic}, \cite{smithies}, and \cite{bg}), I do hope to counteract a false impression students get from current textbooks that Cauchy epitomizes precision and rigor. To me, browsing Cauchy's sprawling oeuvre is like exploring the nooks and crannies of a hyperactive child's tree house, a convoluted structure improvised from scrap lumber and bent nails, remodeled and elaborated over many years. The continuing attraction of the edifice consists in the ingenuity of the creation, the lofty location, and the expansive views from the windows.

\section{First aid.}
Some preliminary observations about the integral~\eqref{eq:quiz} are useful. Cauchy actually includes three positive parameters \(a\), \(b\), and~\(r\) in the integral, thus:
\begin{equation}
 \int_{0}^{\infty} e^{a\cos(bx)} \sin(a\sin(bx)) \frac{x}{x^{2}+r^2} \,dx.
 \label{eq:parameters}
\end{equation}
The integral~\eqref{eq:quiz} is the special case in which \(a=b=r=1\). The number~\(b\) in~\eqref{eq:parameters} can be considered a ``fake parameter'' in the language of \cite{moll}. Indeed, replacing \(bx\) with a new variable~\(u\) produces the equivalent integral
\begin{equation*}
  \int_{0}^{\infty} e^{a\cos(u)} \sin(a\sin(u)) \frac{u}{u^{2}+(br)^2} \,du.
\end{equation*}
Since this integral depends on \(b\) and~\(r\) only through the product~\(br\), there is no loss of generality in setting~\(b\) equal to~\(1\). Once the value of the integral is known when \(b=1\), the general value can be obtained by replacing \(r\) in the answer by~\(br\). 

When \(b=1\), why does the improper integral~\eqref{eq:parameters} converge? Since the expression \( e^{a\cos (x)} \sin(a\sin (x)) \) is an antisymmetric (odd) function of~\(x\), the integral of this quantity over the symmetric interval \([-\pi,\pi]\) is equal to~\(0\). Periodicity implies that the integral of the same expression vanishes over every interval of width \(2\pi\). Therefore 
the magnitude of the integral
\begin{equation*}
 \int_{0}^{R} e^{a\cos (x)} \sin(a\sin (x)) \,dx  
 \end{equation*}
is no more than the magnitude of the integral over half a period, hence
is bounded above by \(\pi e^{a}\) for every~\(R\).
The expression \(x/(x^{2}+r^{2})\) is decreasing when \(x>r\) and has limit~\(0\) when \(x\to\infty\), so Dirichlet's test for integrals \cite[p.~430]{Bromwich} implies that the integral~\eqref{eq:parameters} converges when \(b=1\) (hence for every positive value of the fake parameter~\(b\)). 
 
The convergence is a delicate issue, however, for the improper integral
\begin{equation}
\int_{0}^{\infty} e^{ae^{ibx}} \frac{x}{x^{2}+r^{2}} \,dx,
\label{eq:diverge}
\end{equation}
of which \eqref{eq:parameters} is formally the imaginary part,
actually diverges. To see why, let \(f(x)\) denote the sum of the series \(\sum_{n=1}^\infty \frac{a^n e^{ibnx}}{n!\,ibn}\), absolutely convergent when \(x\)~is a real number, and observe that the derivative \(f'(x)\) equals \(e^{ae^{ibx}} -1\). Adding and subtracting~\(1\) in the integrand, then making use of an explicit antiderivative of \(x/(x^{2}+r^{2})\), and finally integrating by parts shows that
\begin{multline*}
\int_{0}^{R} e^{ae^{ibx}} \frac{x}{x^{2}+r^{2}} \,dx =
\int_{0}^{R} \left(e^{ae^{ibx}}-1\right) \frac{x}{x^{2}+r^{2}} \,dx + \frac{1}{2}\log\left(1+\frac{R^{2}}{r^{2}}\right)  \\
= f(R) \frac{R}{R^{2}+r^{2}} 
- \int_{0}^{R} f(x) \frac{r^{2}-x^{2}}{(x^{2}+r^{2})^{2}}\,dx +  \frac{1}{2}\log\left(1+\frac{R^{2}}{r^{2}}\right).
\end{multline*}
The function~\(f\) is uniformly bounded on the real line, since
\begin{equation*}
\sum_{n=1}^{\infty} \left| \frac{a^n e^{ibnx}}{n!\,ibn} \right| \le \frac{1}{b} \sum_{n=1}^{\infty} \frac{a^{n}}{n!} = \frac{e^{a}-1}{b}.
\end{equation*}
Therefore \(f(R)R/(R^{2}+r^{2})\to 0\) when \(R\to\infty\), and
\begin{equation*}
\lim_{R\to\infty}\left[\int_0^R e^{ae^{ibx}} \frac{x}{x^2+r^2}\,dx  -\frac{1}{2}\log\left(1+\frac{R^2}{r^2}\right)\right] = \int_{0}^{\infty} f(x) \frac{x^{2}-r^{2}}{(x^{2}+r^{2})^{2}}\,dx.
\end{equation*}
Since \(f\)~is bounded, the improper integral on the right-hand side converges absolutely by comparison with the convergent integral \(\int_{0}^{\infty} (x^{2}+r^{2})^{-1}\,dx\). Thus the integral~\eqref{eq:diverge} diverges at a logarithmic rate.

Accordingly, the idea of computing Cauchy's integral~\eqref{eq:parameters} by first evaluating~\eqref{eq:diverge} seems to be a nonstarter. Following Cauchy's lead, I will show nonetheless that the residue method succeeds when~\eqref{eq:parameters} is sneakily realized as the imaginary part of the convergent integral
\begin{equation}
  \int_{-\infty}^{\infty} \frac{1}{2} \left(  e^{ae^{ibx}}-1\right) \frac{x}{x^{2}+r^{2}} \,dx.
\label{eq:sneaky}
\end{equation}

\section{Cauchy's rectangular mallet.}
In 1814, the twenty-five-year-old Cauchy must have been pleased when the French Academy of Sciences accepted his long submission about the evaluation of real definite integrals, especially since his bid for election to that august body had failed the year before~\cite[p.~206]{crosland}. But the actual printing of his article was delayed until 1827, by which time Cauchy had published improved accounts of his theory superseding the first paper.  

I will focus on the following statement from Cauchy's 1814 article \cite[Th\'eor\`eme~1, p.~713]{delay}, paraphrased in modern language:
If \(f(x+yi)\) is holomorphic (that is, complex-analytic) except for some simple poles, and if \(\lim_{y\to\infty} f(x+yi) =0\) for every~\(x\) and \(\lim_{x\to\pm\infty} f(x+yi)=0\) for every~\(y\), then the integral \(\int_{-\infty}^{\infty} f(x)\,dx\) along the real axis equals \(2\pi i\) times the sum of the residues of the function~\(f\) in the upper half-plane. 
Essentially the same statement appears in the 1823 write-up of his calculus lectures \cite[Le\c{c}on 34, p.~136]{resume}. Cauchy's hypotheses actually are not sufficient to guarantee validity of the conclusion, as his mistaken initial evaluation of~\eqref{eq:quiz} reveals. His 1826 correction strengthens one hypothesis~\cite[Th\'eor\`eme VI]{correction} to the still inadequate assumption that
\( (x+yi)f(x+yi)\) tends to~\(0\) when \(y\to\infty\).

But Cauchy's new method does give the right answer for many examples, including this one \cite[p.~758]{delay}:
\begin{equation}
\int_{0}^{\infty} \frac{x\sin (bx)}{x^{2}+1} \,dx = \frac{\pi}{2} e^{-b}
\qquad \text{when \(b>0\)}.
\label{eq:discont}
\end{equation}
Cauchy knows that this formula is correct, since he is aware of an alternative derivation~\cite[p.~100]{laplace} by Pierre-Simon Laplace (1749--1827).
Further on, Cauchy points out \cite[p.~789]{delay} that this integral depends discontinuously on the parameter~\(b\), since the left-hand side of~\eqref{eq:discont} evidently vanishes when \(b=0\), yet the right-hand side reduces to the limiting value~\(\pi/2\). This example challenges a belief cherished by many calculus students that discontinuities appear only in artificial, esoteric situations. 

\begin{figure}
\begin{center}
\small
\begin{tikzpicture}
\draw (-0.5,0) -- (5.5,0);
\draw (0,-0.5) -- (0,3.5);
\draw[thick,decoration={markings, mark=at position 0.5 with {\arrow{>}}}, postaction={decorate}] (1,1) -- (5,1);
\draw[thick,decoration={markings, mark=at position 0.5 with {\arrow{>}}}, postaction={decorate}] (5,1) -- (5,3);
\draw[thick,decoration={markings, mark=at position 0.5 with {\arrow{>>}}}, postaction={decorate}] (1,1) -- (1,3);
\draw[thick,decoration={markings, mark=at position 0.5 with {\arrow{>>}}}, postaction={decorate}] (1,3) -- (5,3);
\fill (1,1) circle (1.5pt);
\fill (5,3) circle (1.5pt);
\draw (1,0.1) -- (1,-0.1) node[anchor=north]{$x_{1}$};
\draw (5,0.1) -- (5,-0.1) node[anchor=north]{$x_{2}$};
\draw (0.1,1) -- (-0.1,1) node[anchor=east]{$y_{1}$};
\draw (0.1,3) -- (-0.1,3) node[anchor=east]{$y_{2}$};
\end{tikzpicture}
\end{center}
\caption{Two paths from \( (x_{1}, y_{1})\) to \((x_{2},y_{2})\).}
\label{fig:paths}
\end{figure}
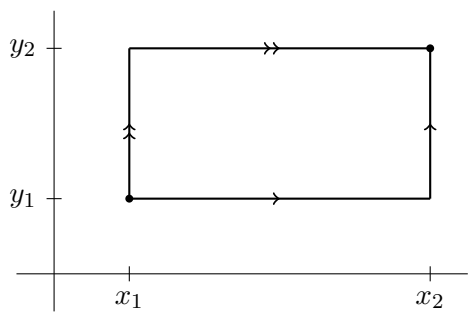

Cauchy's main observation is that computing the integral of a function from a point \( (x_{1},y_{1})\) to a point \( (x_{2},y_{2})\) in two different ways---either along a horizontal path from \( (x_{1},y_{1})\) to \((x_{2},y_{1})\) and a vertical path from \((x_{2},y_{1})\) to \((x_{2},y_{2})\) or, alternatively, along a vertical path from \( (x_{1},y_{1})\) to \((x_{1},y_{2})\) and a horizontal path from \( (x_{1},y_{2})\) to \( (x_{2},y_{2})\)---produces identical answers if the function is holomorphic in the rectangle bounded by the indicated line segments (see Figure~\ref{fig:paths}); and if the function has singularities inside the rectangle, then the two integrals differ by a correction term equal to \(2\pi i\) times the sum of the residues of the function inside the rectangle. My paraphrase is anachronistic: Cauchy does not have the present terminology of line integrals (path integrals), and only in 1826 does he introduce the word ``residue''~\cite{residue}.

\begin{figure}
\begin{center}
\small
\begin{tikzpicture}
\draw (-0.5,0) -- (5.5,0);
\draw (0,-0.5) -- (0,3.5);
\filldraw[fill=gray!20] (1,1) rectangle (5,3);
\begin{scope}[thick,decoration={markings, mark=at position 0.5 with {\arrow{>}}}]
\draw[postaction={decorate}] (1,1.5) -- (5,1.5);
\draw[postaction={decorate}] (1,2) -- (5,2);
\draw[postaction={decorate}] (1,2.5) -- (5,2.5);
\end{scope}
\draw (1,0.1) -- (1,-0.1) node[anchor=north]{$x_{1}$};
\draw (5,0.1) -- (5,-0.1) node[anchor=north]{$x_{2}$};
\draw (0.1,1) -- (-0.1,1) node[anchor=east]{$y_{1}$};
\draw (0.1,3) -- (-0.1,3) node[anchor=east]{$y_{2}$};
\end{tikzpicture}
\end{center}
\caption{Evaluating \( \iint \frac{\partial f}{\partial x}\).}
\label{fig:double}
\end{figure}
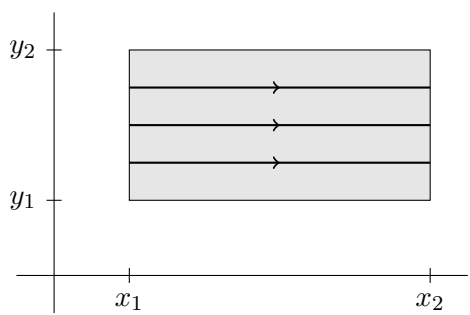

Expressed in today's language, the proof is straightforward. (Students and teachers of multivariable calculus should recognize the argument from the proof of Green's theorem in the plane.\footnote{Although the remarkable George Green (1793--1841) privately published his essay concerning the three-dimensional theorem in 1828, his work became known only after his death. But that is another story \cite{green, gg}.})
To say that a complex-valued function~\(f\) of the two real variables \(x\) and~\(y\) is holomorphic means intuitively that the function depends only on the combination \(x+yi\)
(the complex variable~\(z\)); a more precise statement is that the partial derivative of~\(f\) with respect to~\(y\) equals \(i\)~times the partial derivative with respect to~\(x\). 
Assuming (as Cauchy does implicitly) that the partial derivative \(\partial f/\partial x\) is continuous, apply the fundamental theorem of calculus to rewrite the two-dimensional integral of \(\partial f/\partial x\) over the rectangle as 
\begin{equation*}
\int_{y_{1}}^{y_{2}} f(x_{2},y)-f(x_{1},y) \,dy
\end{equation*}
(see Figure~\ref{fig:double}). 
Similarly, the two-dimensional integral of \(\partial f/\partial y\) over the rectangle equals
\begin{equation*}
\int_{x_{1}}^{x_{2}} f(x,y_{2})-f(x,y_{1}) \,dx.
\end{equation*}
When \(f\) is holomorphic, the second two-dimensional integral equals \(i\)~times the first. The geometric interpretation is that the line integral \( \int f(x,y)\,d(x+yi)\) around the oriented boundary of the rectangle equals~\(0\). 

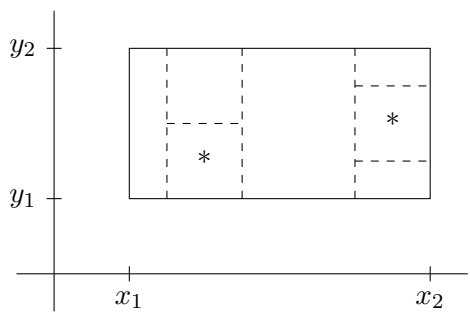
\begin{figure}
\begin{center}
\small
\begin{tikzpicture}
\draw (-0.5,0) -- (5.5,0);
\draw (0,-0.5) -- (0,3.5);
\draw (1,0.1) -- (1,-0.1) node[anchor=north]{$x_{1}$};
\draw (5,0.1) -- (5,-0.1) node[anchor=north]{$x_{2}$};
\draw (0.1,1) -- (-0.1,1) node[anchor=east]{$y_{1}$};
\draw (0.1,3) -- (-0.1,3) node[anchor=east]{$y_{2}$};
\draw (1,1) rectangle (5,3);
\draw (2,1.5) node {*};
\draw (4.5,2) node {*};
\draw[dashed] (1.5,1) -- (1.5,3);
\draw[dashed] (2.5,1) -- (2.5,3);
\draw[dashed] (1.5,2) -- (2.5,2);
\draw[dashed] (4,1) -- (4,3);
\draw[dashed] (4,1.5) -- (5,1.5);
\draw[dashed] (4,2.5) -- (5,2.5);
\end{tikzpicture}
\end{center}
\caption{Localizing the singularities.}
\label{fig:local}
\end{figure}

If \(f\)~has some simple poles (first-order singularities), then a correction term needs to be computed. 
To say that \(f\)~has a simple pole at~\(z_{0}\) with residue equal to the complex number~\(c\) means that the difference \(f(z)-c(z-z_{0})^{-1}\) is holomorphic near~\(z_{0}\) or (equivalently) can be expanded near~\(z_{0}\) in a Taylor series in powers of \((z-z_{0})\). 
Adding and subtracting integrals along suitable line segments and then discarding vanishing integrals over rectangles that avoid the singularities, as indicated in Figure~\ref{fig:local}, reduces the problem to calculating \(\int 1/(z-z_{0}) \,dz\) around a square centered at~\(z_{0}\). 
Making a translation and a dilation converts the problem into showing that \(2\pi i\) equals the value of the counterclockwise line integral
\begin{equation*}
\int \frac{1}{x+yi} \, d(x+yi), \qquad \text{equivalently} \qquad
\int \frac{x-yi}{x^{2}+y^{2}} \,d(x+yi),
\end{equation*}
around the square with vertices at \((\pm 1,\pm 1)\). Symmetry considerations show that the preceding integral equals
\begin{equation*}
4\int_{-1}^{1} \frac{i}{t^{2}+1} \,dt, \quad \text{or} \quad
4i\left[ \arctan(1) - \arctan(-1)\right],
\quad \text{or}
\quad 2\pi i,
\end{equation*}
as claimed. 
(Higher-order singularities can be handled too, as Cauchy makes explicit in~\cite[Th\'eor\`eme~II, p.~131]{correction}, but this refinement is not needed for the main examples.) 

Cauchy's application is to put the bottom edge of the rectangle on the real axis and to let the top and the sides zoom off to infinity. He supposes---wrongly---that if the function vanishes at infinity, then so do the limits of the integrals over the line segments,
and he deduces that the integral of the function over the real axis equals \(2\pi i\) times the sum of the residues of the function in the upper half-plane. 
At this stage in the development of his theory, Cauchy is not thinking about integrals over general simple closed curves: rectangles suffice for evaluating integrals over the real axis.

Issues about domains do not concern Cauchy, for most of his examples involve concrete elementary functions. I will suppose that all functions in question are holomorphic on an open neighborhood of the closed upper half-plane except for finitely many isolated singularities. Using his new invention of singular integrals, Cauchy can allow poles on the real axis, 
but I will assume for simplicity that all the singular points have nonzero imaginary part. With these conventions in force, Cauchy's method can be formalized rigorously as follows.

\begin{theorem}[after Cauchy]
\label{thm:cauchy}
If the improper integral\/ \(\int_{0}^{\infty} f(x+yi) \,dy \) tends to~\(0\) when \(x\to\pm\infty\), and if for each bounded interval~\(I\) the integral\/ \( \int_{I} f(x+yi)\,dx\) tends to~\(0\) when \(y\to\infty\), then\/ \(\int_{-\infty}^{\infty} f(x)\,dx\) equals \(2\pi i\) times the sum of the residues of\/~\(f\) in the upper half-plane.
\end{theorem}

For the proof, fix a small positive number~\(\varepsilon\), and invoke the hypotheses to say that for all sufficiently large positive numbers \(A\) and~\(B\), the integrals \(\int_{0}^{\infty} f(-A+yi) \,dy\) and \(\int_{0}^{\infty} f(B+yi) \,dy\) have absolute value less than~\(\varepsilon/3\), and every singular point of~\(f\) has real part between \(-A\) and~\(B\). Fix such numbers \(A\) and~\(B\), and use the meaning of convergence of an improper integral to deduce that for every sufficiently large positive number~\(C\), the integrals \(\int_{0}^{C} f(-A+yi) \,dy\) and \(\int_{0}^{C} f(B+yi) \,dy\) have absolute value less than~\(\varepsilon/3\). By hypothesis, the number~\(C\) can be chosen additionally large enough that the integral 
\(\int_{-A}^{B} f(x+Ci)\,dx \) has absolute value less than~\(\varepsilon/3\). Consequently, the integral of~\(f\) over the rectangle with opposite corners \( (-A,0)\) and \( (B,C)\) not only equals \(2\pi i\) times the sum of the residues of~\(f\) in the upper half-plane but also differs from \(\int_{-A}^{B} f(x)\,dx\) by less than~\(\varepsilon\). Since \(\varepsilon\)~is arbitrary, the doubly improper integral \(\int_{-\infty}^{\infty} f(x)\,dx\) converges and has the required value.

The problem to which Cauchy has no adequate solution is to specify readily verified conditions on the integrand to guarantee that the hypotheses of the theorem hold. The assumption that \(\lim_{y\to\infty} \int_{I} f(x+yi)\,dx=0\) certainly holds when \(f(x+yi)\) tends to~\(0\) uniformly with respect to~\(x\) when \(y\to\infty\). This assumption may be what Cauchy has in mind when he says that \(\lim_{y\to\infty} f(x+yi)=0\), but he lacks the concept of uniform convergence. The property certainly holds when \(f\)~is a rational function whose denominator has higher degree than the numerator, for then \(|f(x+yi)|\) decays at least as fast as a constant times \(1/|x+yi|\), and this expression is uniformly bounded above by \(1/y\) when \(y\to\infty\). 
A trickier issue is to identify simple but widely applicable conditions on~\(f\) to control the integrals on vertical lines. There is no trouble working with a rational function if the degree of the denominator is at least \(2\) larger than the degree of the numerator, for then \(|f(x+yi)|\) decays at least as fast as a constant times \(1/|x+yi|^{2}\), and 
\begin{equation*}
\int_{0}^{\infty} \frac{1}{x^{2}+y^{2}} \,dy = \frac{\pi}{2|x|},
\end{equation*}
which does tend to~\(0\) when \(x\to\pm\infty\). To handle more general functions requires further analysis, a topic that I address next. 

\section{Whacking Jordan.}
At first sight, Theorem~\ref{thm:cauchy} seems insufficient to evaluate the integral in~\eqref{eq:discont}. That the integral lives on only a part of the real axis is no difficulty: by symmetry, the value equals half the integral over the whole axis. A serious issue, however, is that \(\sin [b(x+yi)]\) blows up as \(y\)~tends to infinity. When \(x=0\), for instance, this function becomes \(\sin(byi)\), which equals the unbounded purely imaginary expression \(i\sinh(by)\). As suggested in the introduction, Cauchy's device for overcoming this obstacle is to express the real function \(\sin(bx)\) as the imaginary part of the complex exponential \( e^{ibx}\) and to view the integral as the imaginary part of
\begin{equation*}
\frac{1}{2} \int_{-\infty}^{\infty} \frac{x}{x^{2}+1} e^{ibx}\,dx.
\end{equation*}
If \(z=x+yi\), then \(|e^{ibz}| = e^{-by}\), an expression that tends to~\(0\) uniformly with respect to~\(x\) when \(y\to\infty\). Since the fraction \(z/(z^{2}+1)\) also tends to~\(0\) uniformly when \(y\to\infty\), the second hypothesis in Theorem~\ref{thm:cauchy} holds. To check the hypothesis about integrals on vertical lines, observe that
\begin{equation*}
\left| \frac{z}{z^{2}+1}\right| = \left| \frac{(z-i)+i}{(z-i)(z+i)}\right| \le 
\frac{1}{|z+i|} + \frac{1}{|z-i|\,|z+i|} \le \frac{1}{|x|} + \frac{1}{x^{2}}.
\end{equation*}
Therefore,
\begin{equation*}
\biggl| \int_{0}^{\infty}  \frac{x+yi}{(x+yi)^{2}+1} e^{ib(x+yi)}  \,dy \biggr|
\le \left( \frac{1}{|x|} + \frac{1}{x^{2}}\right) \int_{0}^{\infty} e^{-by}\,dy.
\end{equation*}
The right-hand side tends to~\(0\) when \(x\to\pm\infty\), the integral \(\int_{0}^{\infty} e^{-by}\,dy\) being finite. Thus both hypotheses of Theorem~\ref{thm:cauchy} are satisfied. The residue of
\begin{equation*}
\frac{1}{2}\cdot\frac{z}{z^{2}+1} e^{ibz}
\end{equation*}
at the singular point~\(i\) equals \(e^{-b}/4\), so Theorem~1 validates formula~\eqref{eq:discont}.

\begin{figure}
\begin{center}
\small
\begin{tikzpicture}
\draw (-3.5,0) -- (3.5,0);
\draw (0,-0.3) -- (0,3.3);
\draw (3,0) node[anchor=north] {$R$};
\draw (-3.1,0) node[anchor=north] {$-R$};
\draw (45:3cm) node[anchor=south west] {$C_{R}$};
\draw[thick] (3,0) arc (0:180:3);
\end{tikzpicture}
\end{center}
\caption{Jordan's contour.}
\label{fig:semi}
\end{figure}
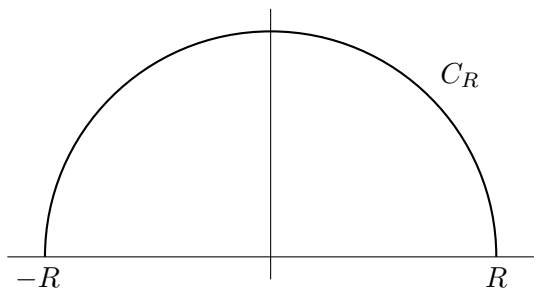

Curiously, a more difficult argument for establishing such integral formulas has been standard in textbooks since the 19th-century \emph{Cours d'analyse} \cite{jordan} of Camille Jordan (1838--1922). Students are taught to evaluate
\begin{equation*}
 \lim_{R\to\infty} \frac{1}{2}  \int_{-R}^{R} \frac{x}{1+x^{2}} e^{ibx} \,dx
\end{equation*}
by closing the contour with a semicircle~\(C_{R}\) in the upper half-plane (see Figure~\ref{fig:semi}). 
This method is trickier than using Cauchy's rectangle. Indeed, I was unable to complete the proof as an undergraduate until a fellow student proposed integrating over a triangle instead of a semicircle. Our instructor was Lars Ahlfors (1907--1996), recipient of the Fields medal at the first award ceremony (in 1936). He indicated that he had not seen triangles used before in this context, so we students published the idea~\cite{eduardo}, naively unaware that rectangles are the original polygonal contrivance from the dawn of the theory. 

A student who attempts to prove that the integral over~\(C_{R}\) has limit equal to~\(0\) when \(R\to\infty\) typically falls short on the first try as follows. Since \(|e^{ib(x+yi)}| = e^{-by}\le 1\), and \(|z/(1+z^{2})|\le R/(R^{2}-1)\) when \(|z|=R\) (as long as \(R>1\)), bounding the integral by the length of the integration path times the maximum of the absolute value of the integrand shows that
\begin{equation}
\biggl| \int_{C_{R}} \frac{z}{1+z^{2}} e^{ibz}\,dz \biggr| \le \pi R \cdot \frac{R}{R^{2}-1}.
\label{eq:CR}
\end{equation}
The indicated upper bound tends to~\(\pi\), not~\(0\), when \(R\to\infty\). This inequality admits an essential improvement, however, for \(e^{-by}\) not only is bounded but actually decays when \(y\)~grows. Since \(y\)~is large only on part of the semicircle~\(C_{R}\), care is needed to exploit this extra information.

The second try is to parametrize the semicircle by setting~\(z\) equal to \(Re^{i\theta}\). Moving the absolute value signs inside the integral yields the following upper bound for the left-hand side of~\eqref{eq:CR}:
\begin{equation}
\frac{R^{2}}{R^{2}-1} \int_{0}^{\pi} e^{-bR\sin(\theta)}\,d\theta. 
\label{eq:nonuniform}
\end{equation}
The factor in front of the integral tends to~\(1\) when \(R\to\infty\), so what needs to be shown is that the integral with respect to~\(\theta\) has limit equal to~\(0\). When \(0<\theta<\pi\), the quantity \(b\sin(\theta)\) is positive, so the integrand \(e^{-bR\sin(\theta)}\) has limit~\(0\). 
But this limit is not uniform with respect to~\(\theta\), a complication that could trip up even professional mathematicians as recently as the early 20th century.
Edmund Taylor Whittaker (1873--1956) waves his hands unconvincingly in the first edition of \emph{A Course of Modern Analysis} \cite[p.~86]{first}, claiming  nonsensically that since the original integrand \(ze^{ibz}/(1+z^{2})\) ``is infinitesimal compared with \(1/z\) at points on~\(\gamma\) [\(=C_{R}\)], the integral round~\(\gamma\) is infinitesimal compared with \(\int_{\gamma} |dz/z|\) or \(2\pi\), and is therefore zero.'' The integrand is certainly not small in comparison to \(1/z\) on the real axis! Nonetheless, the required convergence of the integral~\eqref{eq:nonuniform} 
to~\(0\) when \(R\to\infty\) does follow directly from any one of several propositions that are standard nowadays: the bounded convergence theorem,  the dominated convergence theorem, and the monotone convergence theorem. 

The preceding analysis shows, more generally, that if \(f\)~is a function for which \(|zf(z)|\) is bounded when \(|z|\) is sufficiently large, then \(\int_{C_{R}} f(z) e^{ibz}\,dz\) tends to zero when \(R\to\infty\). In particular, \(f\)~can be any rational function that vanishes at infinity.

This argument can be sharpened by leveraging the fast decay of the exponential to avoid using any quantitative information about the decay rate of the other factor.
Jordan observes \cite[\S289]{jordan}
for a general function~\(f\) that
\begin{equation*}
\biggl| \int_{C_{R}} f(z)e^{iz}\,dz \biggr| \le \sup_{0<\theta<\pi} |f(Re^{i\theta})| 
\int_{0}^{\pi} e^{-R\sin(\theta)} R\,d\theta.
\end{equation*}
(Adjusting the argument to hold for \(e^{ibz}\) in place of \(e^{iz}\) is a simple matter of rescaling the variable.) 
Since \(\sin(\pi-\theta)=\sin(\theta)\), the \(\theta\)~integral equals twice the integral from \(0\) to~\(\pi/2\). Now \(\sin(\theta) \ge (2/\pi)\theta\) when \(0\le \theta\le \pi/2\), so an upper bound for the \(\theta\)~integral is 
\begin{equation*}
2\int_{0}^{\pi/2} e^{-(2/\pi)R\theta} R\,d\theta,
\qquad \text{or} \qquad \pi(1-e^{-R}).
\end{equation*}
Accordingly, if \(f(Re^{i\theta})\to 0\) when \(R\to\infty\), and if the convergence is uniform with respect to the angle~\(\theta\), then \(\int_{C_{R}} f(z)e^{iz}\,dz \) has limit~\(0\) too. 

Rarely needed in practice, this refinement is dubbed  ``Jordan's lemma'' in textbooks. The person responsible for naming the lemma  seems to be George Neville Watson (1886--1965) in his textbook~\cite[\S30]{watson}. Although historically accurate, the nomenclature is unfortunate, for an identically named algebraic proposition exists in the theory of invariants (see, for instance, \cite[Appendix~III]{invariant}). The name of the analytic lemma is recorded in the second edition of \emph{A Course of Modern Analysis}~\cite{second}, presumably due to Watson's collaboration with his teacher on the revision; the wide influence of this book, still in print, has ensured the permanence of the terminology. Jordan himself died at the age of~84 during the interval between the third edition and the fourth edition of 
this treatise, familiarly known as ``Whittaker and Watson.''

Watson's own textbook evaluates Cauchy's integral~\eqref{eq:parameters} by ``arguments similar to those used in proving Jordan's lemma'' \cite[p.~62]{watson}, but the exposition requires some contortions, since the lemma does not apply as formulated. Thus there is a need to upgrade Jordan's lemma. I offer a replacement in the spirit of Cauchy's initial work on definite integrals. Although not explicitly covered by Cauchy's writings of 1825--1826, the statement is in the penumbra. 

The integrand in Jordan's lemma is a product of two functions of different character. Accordingly, I consider functions \(f_{1}\) and~\(f_{2}\) that are holomorphic in a neighborhood of the closed upper half-plane and have finitely many singularities, all located in the open upper half-plane. 

\begin{theorem}[after Cauchy]
Suppose that when \(y\to\infty\), the function \(f_{1}(x+yi)\) tends to~\(0\) uniformly with respect to~\(x\) in an arbitrary bounded interval~\(I\), and the integral \(\int_{I} |f_{2}(x+yi)|\,dx\) stays bounded. Suppose that when \(x\to\pm\infty\), the function \(f_{2}(x+yi)\) tends to~\(0\) uniformly with respect to~\(y\), and the integral \(\int_{0}^{\infty} |f_{1}(x+yi)|\,dy\) stays bounded. Then \(\int_{-\infty}^{\infty} f_{1}(x)f_{2}(x)\,dx\) equals \(2\pi i\) times the sum of the residues of \(f_{1}f_{2}\) in the upper half-plane. 
\label{lemma:two}
\end{theorem}

This statement is an immediate corollary of Theorem~\ref{thm:cauchy}, for the hypotheses imply that the product function~\(f_{1}f_{2}\) satisfies the conditions of that theorem. 
The traditional Jordan lemma is the special case in which \(f_{1}(z)\)~is an exponential function of the form \(e^{ibz}\) (where \(b>0\)) and \(f_{2}(z)\) vanishes when \(z\to\infty\).

There are many other interesting functions to which Theorem~\ref{lemma:two} applies. Indeed, let \(g(w)\) be an arbitrary power series \(\sum_{n=1}^{\infty} c_{n} w^{n}\) that has radius of convergence greater than~\(1\) and lacks a constant term. I claim that the composite function \(g(e^{ibz})\) will serve for~\(f_{1}(z)\) in Theorem~\ref{lemma:two} when \(b>0\). Indeed,
\begin{equation*}
|g(e^{ib(x+yi)})| \le e^{-by} \sum_{n=1}^{\infty} |c_{n}| 
\qquad \text{when \(y>0\)}.
\end{equation*}
Now \( \sum_{n=1}^{\infty} |c_{n}| \) converges, because 
 \(\sum_{n=1}^{\infty} c_{n} = g(1)\), and every power series converges absolutely inside the open disk of convergence. Accordingly, the function \(g(e^{ib(x+yi)})\) tends to~\(0\) uniformly with respect to~\(x\) when \(y\to\infty\). The preceding inequality additionally implies that
\begin{equation*}
\int_{0}^{\infty} |g(e^{ib(x+yi)})| \,dy \le \frac{1}{b} \sum_{n=1}^{\infty} |c_{n}|,
\end{equation*}
the finite upper bound being independent of~\(x\). Thus the function \(g(e^{ibz})\) does satisfy both hypotheses required of \(f_{1}(z)\) in Theorem~\ref{lemma:two}.

In particular, choosing \(e^{aw}-1\) for the function \(g(w)\) shows that \(f_{1}(z)\) can be taken to be \(e^{ae^{ibz}}-1\) in Theorem~\ref{lemma:two}. Letting \(f_{2}(z)\) be \(z/(z^{2}+r^{2})\) reveals that 
the value of the integral~\eqref{eq:sneaky} is \(2\pi i\) times the residue at~\(ir\) of the function
\begin{equation*}
\frac{1}{2} \cdot  \frac{z}{z^{2}+r^{2}} \left( e^{ae^{ibz}} -1 \right),
\qquad \text{the residue being}
\qquad \frac{1}{4} \left( e^{ae^{-br}} -1 \right).
\end{equation*}
Consequently, the integral~\eqref{eq:parameters} equals 
\begin{equation*}
\frac{\pi}{2} \left( e^{ae^{-br}} -1\right);
\qquad \text{and the integral \eqref{eq:quiz} equals}
\qquad \frac{\pi}{2} \left( e^{1/e} -1\right).
\end{equation*}
Cauchy's blunder in his first attempt is the failure to account for the term~\(-1\) needed to make~\(g(0)\) equal to~\(0\). 

\section{Driving the point home.}
Cauchy's marvelous tool for computing definite integrals has remained useful into modern times, notwithstanding the development of automated computation. One example dear to my heart is the following formula from my mother's Ph.D. dissertation on theoretical physics \cite{mother}, directed by Herman Feshbach:
\begin{equation*}
 \int_0^\infty \frac{1}{x(1+x^2)}\log\left| \frac{x+\sqrt{3}}{x-\sqrt{3}}\right|\,dx = \frac{\pi^2}{6}.
\end{equation*}
The integrand has a removable singularity when \(x=0\) and an integrable singularity when \(x=\sqrt{3}\). 
The relatively fast decay of the integrand at infinity yields an easy verification of the main hypotheses of Theorem~\ref{thm:cauchy}, and the indicated result follows after a bit of care to define a suitable branch of the complex logarithm function and a finesse to handle innocuous singularities on the real axis. Yet both Maple and \textit{Mathematica} beg the question by evaluating the integral in terms of the dilogarithm function, which itself is defined as an integral, and neither software program succeeds in simplifying\footnote{An amusing exercise for human readers is to massage the computer's output into the required simple form by applying two known identities for the dilogarithm function \cite[Theorem 2.6.1 and Equation 2.6.6]{special}.} the answer to \(\pi^2/6\).

To emphasize the continuing strength and value of Cauchy's simple rectangle method, I offer the following sampler of additional formulas that can be deduced from Theorem~\ref{lemma:two}. As in the integral~\eqref{eq:parameters}, the parameters \(a\), \(b\), and~\(r\) represent  arbitrary positive numbers. In the final three formulas, the additional positive parameter~\(c\) is assumed to have a value greater than~\(a\).

\begin{align}
\int_{0}^{\infty} \sin(a\cos (bx)) \sinh(a\sin (bx)) \frac{x}{x^{2}+r^{2}} \,dx &= \frac{\pi}{2} \left( 1-\cos\left(\frac{a} {e^{br}}\right)\right)
\label{hyper-one} \\
\int_{0}^{\infty} \cos(a\cos  (bx)) \sinh(a\sin (bx)) \frac{x}{x^{2}+r^{2}}\,dx &=\frac{\pi}{2} \sin\left( \frac{a}{e^{br}}\right)
\label{hyper-two}
\\
\int_{0}^{\infty} \frac{x\sin (bx)}{(x^{2}+r^{2}) (a^2+2ac \cos (bx) + c^{2})}\,dx &= \frac{\pi}{2c} \cdot \frac{1}{a+ce^{br}}
\label{denom-one}
\\
\int_{0}^{\infty}  \frac{\log(a^{2} +2ac\cos(b x)+ c^{2})}{x^{2}+r^{2}} \,dx &= \frac{\pi}{r} \log\left( c +\frac{a} {e^{br}}\right) 
\label{log-one}
\\
 \int_{0}^{\infty} \frac{x}{x^{2}+r^{2}} \log\frac{a^{2} +2ac \sin (bx) + c^{2}}{a^2-2ac\sin(bx)+c^2} \,dx &= 2\pi \arctan\left( \frac{a}{ce^{br}}\right)
 \label{log-two} 
\end{align}

Which of these formulas can you prove? Here is your assessment rubric:
\begin{description}
\item [One correct]
You beat both Maple~18 and \textit{Mathematica}~11, which cannot solve any of these problems. Indeed, the computer programs have to be coached even to produce accurate numerical approximations of these slowly converging integrals.

\item [Two correct]
You are on a par with the tables of Gradshteyn and Ryzhik, which to the best of my knowledge contain only \eqref{denom-one} and~\eqref{log-one} \cite[3.792(13), 4.322(10)]{table}.

\item [Three correct]
You outperform the tables of Bierens de Haan, which contain a correct version of~\eqref{hyper-one} \cite[375(1)]{bierens} 
but erroneous versions of \eqref{hyper-two} and~\eqref{denom-one} \cite[375(3), 192(2)]{bierens}.

\item [Four correct]
You may apply for a job as assistant to Cauchy, whose own supplementary list of integrals that his method handles ``without difficulty'' includes the left-hand sides of \eqref{log-one} and~\eqref{log-two} explicitly and \eqref{hyper-one} and~\eqref{hyper-two} implicitly \cite[pp.~88--89]{Cauchy1825}.

\item [Five correct]
Congratulations! You are in a position now to establish ``an infinity of other [examples]'' \cite[p.~88]{Cauchy1825} by the rectangle method. 
\end{description}

\end{document}